\documentclass[11pt]{amsart}
\input xy
\xyoption{all}
\usepackage{graphicx}
\usepackage{latexsym}
\usepackage{amsmath}
\usepackage{amssymb}

\newtheorem{theorem}{Theorem}[section]

\newtheorem{proposition}[theorem]{Proposition}
\newtheorem{conj}[theorem]{Conjecture}
\newtheorem{cor}[theorem]{Corollary}

\theoremstyle{definition}

\theoremstyle{remark}
\numberwithin{figure}{section}

\newcommand{\zz}{{\mathbb Z}}
\newcommand{\ezt}{E_2({\mathbb Z}[t])}
\newcommand{\slzt}{SL_2({\mathbb Z}[t])}

\newcommand{\slzti}{SL_2({\mathbb Z}[t,t^{-1}])}

\newcommand{\slft}{SL_2(F[t])}

\newcommand{\zp}{{\mathbb P}}
\newcommand{\slfti}{SL_2(F[t,t^{-1}])}
\newcommand{\zq}{{\mathbb Q}}
\newcommand{\fp}{{\mathbb F}_2}
\newcommand{\slfpt}{SL_2(\fp[t])}
\newcommand{\slfpti}{SL_2(\fp[t,t^{-1}])}
\newcommand{\slfp}{SL_2(\fp)}

\begin{document}

\title{Homology and finiteness properties of $\slzti$}

\author{Kevin P.~Knudson}
\address{Department of Mathematics \& Statistics, Mississippi State University, Mississippi State, MS 39762}
\email{knudson@math.msstate.edu}

\date{\today}

\begin{abstract}
We show that the group $H_2(\slzti;\zz)$ is not finitely generated, answering a question mentioned by Bux and Wortman in \cite{bux}.
\end{abstract}

\maketitle

\section{Introduction}  In \cite{bux}, Bux and Wortman give a geometric proof that the group $\slzti$ is not finitely presented (nor is it of type $\text{FP}_2$).  They achieve this by letting the group act on a product of Bruhat--Tits trees.  In their introduction, they point out the following unsolved question:

\begin{enumerate}
\item[] Is $H_2(\slzti;\zz)$ finitely generated?
\end{enumerate}

In this note, we prove that the answer is no.

\medskip

\noindent {\bf Theorem \ref{notfingen}}.  {\em The group $H_2(\slzti;\zz)$ is not finitely generated}.

\medskip

We conjecture that Theorem \ref{notfingen} is true for all $i\ge 1$; that is $H_i(\slzti;\zz)$ is not finitely generated for all $i\ge 1$.  Our method of proof suggests how one might prove this for all $i\ge 2$, but the method breaks down for $i=1$.  As a consequence, we would have that $\slzti$ is not finitely generated, an open question also mentioned in \cite{bux}.  For $n\ge 3$, however, Suslin \cite{suslin} proved that the group $SL_n(\zz[t,t^{-1}])$ is finitely generated by elementary matrices.  The question of whether $\slzti$ is generated by elementary matrices is an open question.

We note also that in a recent preprint \cite{abramenko}, Abramenko has proved various interesting results about non-generation of $SL_2(R)$ by elementary matrices for a broad class of rings $R$.  The techniques he uses are similar and more general than what we use below, but do not apply to the ring $R=\zz[t,t^{-1}]$.

The main idea behind the proofs of this result is the following.  Let $F$ be a field.  The group $\slfti$ admits an amalgamated free product decomposition $$\slfti \cong \slft \ast_{\Gamma(F)} \slft,$$ where $$\Gamma(F)=\biggl\{ \left(\begin{array}{rr}
                                                                                            a & b \\
                                                                                            tc & d
                                                                                            \end{array}\right) : a,b,c,d \in F[t]\biggr\}$$
(see e.g. \cite{knudson1}).  Consider the case $F=\zq$.  Letting $P=\slzt \cap \Gamma(\zq)$, we show that the group $\slzti$ admits an amalgamated free product decomposition
$$\slzti=\slzt \ast_{P} \slzt.$$

A portion of the homology of $\slzt$ was calculated by the author in \cite{knudson2}.  In particular, it is not finitely generated in any positive degree.  Similarly, the homology of $\slfpt$ is not finitely generated in any positive degree, and a consideration of the commutative diagram
of Mayer--Vietoris sequence associated to the amalgamated free product decompositions of $\slzti$ and $\slfpti$ yields Theorem \ref{notfingen}.

\medskip

\noindent {\em Acknowledgments}.  The author thanks Kai-Uwe Bux and Kevin Wortman for valuable comments about this manuscript.  Thanks are due also to Peter Abramenko for finding an error in a previous version and to Roger Alperin for pointing out \cite{krstic}.

\section{Amalgamated Free Products}\label{elementary}

Let $F$ be a field and consider the commutative diagram
$$\xymatrix{
\Gamma(F) \ar[r]^{j_1}\ar[d]_{j_2} & \slft \ar[d]^{i_1} \\
\slft \ar[r]^{i_2} & \slfti}$$
where $\Gamma(F)$ consists of those matrices in $\slft$ that are upper triangular modulo $t$, $j_1$ and $j_2$ are the standard inclusion maps, $i_1$ is the inclusion, and $i_2:\slft\to \slfti$ is defined by
$$i_2(A) = \left(\begin{array}{cc}
                   1 & 0  \\
                   0 & t
                   \end{array}\right) A \left(\begin{array}{cc}
                                                1 & 0 \\
                                                0 & t^{-1}
                                                \end{array}\right).$$
Using standard arguments for groups acting on trees (see e.g. \cite{serre}), one can prove that $\slfti$ admits an amalgamated free product decomposition
$$\slfti \cong \slft \ast_{\Gamma(F)} \slft.$$

Now consider the case $F=\zq$ and let $P=\slzt\cap \Gamma(\zq)$. 

\begin{proposition}\label{amalgam} There is an amalgamated free product decomposition
$$\slzti \cong \slzt \ast_P \slzt.$$
\end{proposition}

\begin{proof} Let $H$ be the subgroup of $\slzti$ generated by $\slzt$ and $i_2(\slzt)$.  Note that $$P=\slzt\cap \Gamma(\zq) = i_2(\slzt)\cap \Gamma(\zq).$$ Indeed, $\Gamma(\zq)=SL_2(\zq[t])\cap i_2(SL_2(\zq[t]))$, so that
\begin{eqnarray*}
\slzt \cap \Gamma(\zq) & = & \slzt\cap SL_2(\zq[t])\cap i_2(SL_2(\zq[t])) \\
                       & = & \slzt\cap i_2(SL_2(\zq[t])) \\
                       & = & \slzt\cap i_2(\slzt) \\
                       & = & i_2(\slzt)\cap \Gamma(\zq).
\end{eqnarray*}
  By \cite{serre}, p.~6, the evident map
$$\slzt \ast_P \slzt \to H$$ is an isomorphism.  It remains to show that $H=\slzti$, but this follows easily from Theorem 3.4 of \cite{abramenko} by taking $R=\zz[t]$, $\pi=t$.
\end{proof}

\noindent {\em Remark.}  Proposition \ref{amalgam} is true for any PID $R$, again by using Theorem 3.4 of \cite{abramenko}.

\section{Homology}\label{homology}
In this section, we prove the following result.

\begin{theorem}\label{notfingen} The group $H_2(\slzti;\zz)$ is not finitely generated.
\end{theorem}

\begin{cor}\label{notfinpres} The group $\slzti$ is not finitely presented.
\end{cor}

\begin{proof} This follows from the fact that $H_2(\slzti;\zz)$ is not finitely generated.
\end{proof}

\noindent {\em Remark.} Corollary \ref{notfinpres} was proved by Bux and Wortman \cite{bux} using group actions on trees and also by Krsti\'c and McCool \cite{krstic} by more elementary group theoretic means.  In fact, any subgroup of $GL_2(\zz[t,t^{-1}])$ containing the subgroup of elementary matrices is not finitely presented.

\medskip

\noindent {\em Proof of Theorem \ref{notfingen}}.  Consider the map
$$\pi:\slzti\to\slfpti$$
induced by the map $\zz[t,t^{-1}]\to \fp[t,t^{-1}]$ which reduces coefficients modulo $2$.  To prove the result, it suffices to show that the image of $$\rho_*:H_2(\slzti;\zz)\to H_2(\slfpti;\zz)$$ is not finitely generated.

Consider the commutative diagram of Mayer--Vietoris sequences
$$\xymatrix{
H_2(P;\zz)\ar[r]\ar[d] & H_2(\slzt;\zz)\oplus H_2(\slzt;\zz)\ar[r]\ar[d] & H_2(\slzti;\zz)\ar[d] \\
H_2(\Gamma(\fp);\zz)\ar[r] & H_2(\slfpt;\zz)\oplus H_2(\slfpt;\zz)\ar[r] & H_2(\slfpti;\zz)}$$
The author provided a partial calculation of the homology of $\slzt$ and $\slfpt$ in \cite{knudson2}.  We must refine those, however, to prove the result.

For a ring $R$ denote by $B(R)$ the group of invertible upper-triangular matrices over $R$.  We have an amalgamated free product decomposition
$$\ezt \cong SL_2(\zz) \ast_{B(\zz)} B(\zz[t])$$ which yields a short exact sequence, for all $i\ge 1$,
$$0\to H_i(B(\zz)) \to H_i(SL_2(\zz))\oplus H_i(B(\zz[t])) \to H_i(\ezt)\to 0.$$  The homology of $B(\zz)$ is easily calculated using the Lyndon--Hochschild--Serre (LHS) spectral sequence associated to the extension
$$0\to \zz\to B(\zz)\to \zz^\times \to 1;$$ the result is
$$H_i(B(\zz)) \cong \begin{cases}
                         \zz & i=0 \\
                         \zz\oplus\zz_2 & i=1 \\
                         \zz_2 & i\ge 2.
                         \end{cases}$$

The homology of $B(\zz[t])$ is a bit more complicated.  Again, we consider the LHS spectral sequence associated to the extension
$$0\to t\zz[t] \to B(\zz[t]) \stackrel{t=0}{\to} B(\zz) \to 1.$$  This sequence has $E^2$-term equal to
$$E^2_{p,q}= H_p(B(\zz),H_q(t\zz[t])) = H_p(B(\zz))\otimes H_q(t\zz[t]),$$ the latter equality following from the fact that $B(\zz)$ acts trivially on $t\zz[t]$.  The lower-left corner of this looks as follows
$$\begin{array}{|ccc}
\bigwedge\nolimits^2 t\zz[t] &                   &                                \\
t\zz[t]    & t\zz[t]\oplus \zz_2\otimes t\zz[t] & \zz_2\otimes t\zz[t]  \\
\zz       &   \zz\oplus\zz_2     &   \zz_2    \\ \hline
\end{array}$$
Since the group $E^2_{2,1}=\zz_2\otimes t\zz[t]$ is torsion, and $E^2_{0,2} = \bigwedge^2 t\zz[t]$ is torsion-free, the differential $d^2:E^2_{2,1}\to E^2_{0,2}$ is the zero map, and since all differentials coming from the bottom row vanish (since the extension is split), the lower-left corner of the $E^2$-term is that of $E^\infty$.  It follows that $$H_2(B(\zz[t])) \cong \bigwedge\nolimits^2 t\zz[t] \oplus t\zz[t]\oplus t\fp[t] \oplus \zz_2.$$
Putting this into the short above exact sequence above yields
$$H_2(\ezt)\cong \bigwedge\nolimits^2 t\zz[t] \oplus t\zz[t]\oplus t\fp[t].$$

Now consider the same calculation for $\slfpt\cong \slfp \ast_{B(\fp)} B(\fp[t])$.  Since $\fp^\times$ is the trivial group, we have $B(\fp)\cong \fp$.  It follows that $H_i(B(\fp)) = H_i(\fp)$, $i\ge 0$.  For ease of computing the map on homology induced by the map $\ezt\to\slfpt$, we use the LHS spectral sequence associated to the extension
$$0\to t\fp[t]\to B(\fp[t]) \stackrel{t=0}{\to} B(\fp)\to 1$$ to compute homology.  We find that
$$H_2(B(\fp[t])) \cong \bigwedge\nolimits^2 t\fp[t] \oplus \fp\otimes t\fp[t].$$  Putting this into the short exact sequence, we obtain
$$H_2(\slfpt) \cong \bigwedge\nolimits^2 t\fp[t] \oplus \fp\otimes t\fp[t].$$

\medskip

\noindent {\em Remark.}  Note that the preceding calculations are really just spectral sequence presentations of the K\"unneth formula.  Indeed, since the extensions used are split and the actions of the quotients on the kernels are trivial, the groups $B(\zz)$ and $B(\zz[t])$ are {\em direct} products of the factors involved. This implies the direct sum decompositions for the $H_2$-groups, which would not follow automatically otherwise.  We have presented it this way, however, because of the necessary use of spectral sequences which follows.

\medskip

The naturality of all the spectral sequences involved (or, equivalently, the naturality of the K\"unneth formula) allows us to deduce that the homology homomorphism induced by $\ezt\to \slfpt$ is as follows:
$$\xymatrix{
H_2(\ezt) & = & t\zz[t]\ar[d]_{\text{mod}\,2} & \oplus & \bigwedge\nolimits^2 t\zz[t]\ar[d]_{\text{mod}\,2} & \oplus & t\fp[t]\ar[d]  \\
H_2(\slfpt) & = & t\fp[t] & \oplus & \bigwedge\nolimits^2 t\fp[t] & \oplus & 0}$$  Since this map factors through $H_2(\slzt)$, we conclude that $H_2(\slzt)$ is not finitely generated.
Denote by $u_i$ the image of the element $t^i\in t\zz[t] \subset H_2(\ezt)$ in $\slzt$, and by $v_{jk}$ the image of the element $t^j\wedge t^k\in \bigwedge\nolimits^2 t\zz[t]$.  Then the $u_i$ and $v_{jk}$ generate an infinitely generated subgroup of $H_2(\slzt)$.

Let us now turn our attention to the computation of the homology of $\slfpti$.  To this end, we must compute the homology of the group $\Gamma(\fp)$ consisting of those matrices in $\slfpt$ which are upper-triangular modulo $t$.  We have the split extension
$$1\to K(\fp) \to \Gamma(\fp) \stackrel{t=0}\to B(\fp)\to 1,$$ where $K(\fp)$ consists of those matrices which are congruent to the identity modulo $t$.  Let $C$ be the subgroup of $K$ defined by
$$C=\Biggl\{\left(\begin{array}{cc} 1 & tp(t) \\ 0 & 1 \end{array}\right): p(t)\in F[t]\Biggr\}.$$  For each $s\in\zp^1(F)$, set
$$m(s) = \begin{cases}
            \left(\begin{array}{cc} 1 & 0 \\ s & 1 \end{array}\right) & s\ne\infty \\
            \left(\begin{array}{cc} 0 & 1 \\ 1 & 0 \end{array}\right) & s=\infty.
            \end{cases}$$
Then we have a free product decomposition (see, e.g. \cite{knudson2}) $$K= \ast_{s\in\zp^1(F)} m(s)Cm(s)^{-1}.$$  In the case of $F=\fp$, there are only three factors, corresponding to $0,1,\infty$, and the homology of $K$ decomposes as a direct sum of three copies of that of $C$, say
$$H_k(K) = M_0^k \oplus M_1^k \oplus M_\infty^k,$$ where each $M_i^k = H_k(C)$.  The LHS spectral sequence for computing the homology of $\Gamma(\fp)$ has $E^2$-term
$$E^2_{jk}= H_j(B(\fp),H_k(K)).$$  To compute these terms, suppose $k\ge 1$ and note that we have another spectral sequence with $E^2$-term
$$E^2_{pq}=H_p(\fp^\times,H_q(\fp,H_k(K))) \Rightarrow H_{p+q}(B(\fp),H_k(K)).$$  If $p>0$, these terms vanish since $\fp^\times$ is the trivial group.  When $p=0$ we obtain
$$H_0(\fp^\times,H_q(\fp,H_k(K))) = H_q(\fp,H_k(K)) = \begin{cases}
                                                           H_q(\fp,M_0^k) & q>0 \\
                                                           M_0^k \oplus M_\infty^k & q=0,
                                                           \end{cases}$$
the latter equality following from Shapiro's Lemma, since the action of $\fp$ on $\zp^1(\fp)-\{0\}$ is transitive and the stabilizer of $\infty$ is trivial, and $\fp$ acts trivially on $M_0^k$ (for more details, see \cite{knudson1}).  From this we deduce that
$$H_i(B(\fp),H_k(K)) = \begin{cases}
                           M_0^k\oplus M_\infty^k & i=0  \\
                           H_i(\fp)\otimes M_0^k & i>0.
                           \end{cases}$$

Now, to calculate the low-dimensional homology of $\Gamma(\fp)$, we return to the spectral sequence:
$$\begin{array}{|ccc}
M_0^2 \oplus M_\infty^2 &  \fp\otimes M_0^2  &                                \\
M_0^1 \oplus M_\infty^1 & \fp\otimes M_0^1 &   0 \\
\zz       &   \fp     &   0    \\ \hline
\end{array}$$
Since the differentials starting on the bottom row vanish, and no other differentials affect $E^2_{0,2}$ and $E^2_{1,1}$, we have the following
\begin{eqnarray*}
H_2(\Gamma(\fp)) & = & 0 \oplus \fp\otimes M_0^1 \oplus M_0^2\oplus M_\infty^2  \\
                 & = & 0 \oplus t\fp[t] \oplus \bigwedge\nolimits^2 t\fp[t]\oplus \bigwedge\nolimits^2 t\fp[t].
\end{eqnarray*}

It remains to calculate the map $H_2(\Gamma(\fp))\to H_2(\slfpt)\oplus H_2(\slfpt)$ in the Mayer--Vietoris sequence for computing $H_2(\slfpti)$.
The domain consists of elements of the form $(x,p\wedge q,r\wedge s)\in t\fp[t] \oplus \bigwedge\nolimits^2 t\fp[t]\oplus \bigwedge\nolimits^2 t\fp[t]$.  The target consists of two copies of $H_2(SL_2(\fp))\oplus t\fp[t] \oplus \bigwedge\nolimits^2 t\fp[t]$.  Of course, one of these copies is twisted by the isomorphism $i_2$, and $H_2(SL_2(\fp))=0$.  Taking this into account the map in question is
$$(x,p\wedge q,r\wedge s)\mapsto (x,p\wedge q + r\wedge s, x, p\wedge q + tr\wedge ts).$$  It follows that $H_2(\slfpti)$ is not finitely generated (a fact proved by Stuhler \cite{stuhler}), and that the image of the map
$H_2(\slfpt)\oplus H_2(\slfpt) \to H_2(\slfpti)$ contains a copy of $$t\fp[t]\oplus t\fp[t]/\{(x,x):x\in t\fp[t]\},$$ which is not finitely generated.

Finally, consider the commutative diagram
$$\xymatrix{
H_2(\slzt)\oplus H_2(\slzt)\ar[r]\ar[d] & H_2(\slzti)\ar[d]^{\rho_*} \\
H_2(\slfpt)\oplus H_2(\slfpt)\ar[r] & H_2(\slfpti)}$$
We have the elements $u_i\in H_2(\slzt)$ which map to $t^i\in t\fp[t]\subset H_2(\slfpt)$.    Let $\alpha_i$ be the image of $u_i$ in $H_2(\slzti)$ under the map induced by the inclusion $i_1:\slzt\to H$.  Then the $\alpha_i$  map to a non-finitely generated subgroup of $H_2(\slfpti)$.  It follows that the group $H_2(\slzti)$ is not finitely generated. \hfill $\qed$

\section{Remarks}\label{remarks}

The proof of Theorem \ref{notfingen} suggests the following conjecture.

\begin{conj} For all $i\ge 1$, $H_i(\slzti)$ is not finitely generated.
\end{conj}

Note that, with careful bookkeeping, it should be possible to compute $H_i(\slfpti)$ for all $i\ge 1$.  The spectral sequences for computing $H_i(\slfpt)$ are rather complicated, though, and one suspects it would be a tedious calculation.  The same is true for $\ezt$. Still, it is evident that all the pieces are not finitely generated and so one suspects the same is true for $H_i(\slzti)$ and $H_i(\slfpti)$.

The argument breaks down for $i=1$, however.  We have a commutative diagram
$$\xymatrix{
H_1(P)\ar[r]\ar[d] & H_1(\slzt)\oplus H_1(\slzt)\ar[r]\ar[d] & H_1(\slzti)\ar[r]\ar[d] & 0 \\
H_1(\Gamma(\fp))\ar[r] & H_1(\slfpt)\oplus H_1(\slfpt)\ar[r] & H_1(\slfpti)\ar[r] & 0}$$
The bottom row is
$$0\to \fp[t]\oplus t\fp[t] \to H_1(SL_2(\fp))\oplus \fp[t]\oplus H_1(SL_2(\fp))\oplus \fp[t] \to H_1(\slfpti)\to 0$$ and we see that $H_1(\slfpti)$ is finitely generated (indeed, the group $\slfpti$ is finitely generated).  In \cite{knudson2}, the author showed that $H_1(\slzt)$ contains an infinite rank free abelian group orthogonal to the image of $H_1(\ezt)$ (which is also not finitely generated).  Also, the group $P$ fits into a sequence $$1\to U\to P\stackrel{t=0}{\to} B(\zz)\to 1$$ with $U\subset K=\ast_{s\in\zp^1(\zq)} m(s)Cm(s)^{-1}$.  The Kuro\v{s} Subgroup Theorem \cite{kuros} implies that $U$ is also a free product of subgroups of the $m(s)Cm(s)^{-1}$ and a free group.  This calculation seems hopelessly complicated, however, and so we leave the question of the finite generation of $H_1(\slzti)$ unanswered.

\end{document}